\documentclass[11pt,fullpage,doublespace]{article}
\usepackage{cite}
\usepackage[affil-it]{authblk}
\usepackage{amssymb}
\usepackage[left=1in,top=1in,right=1in,bottom=1in,nohead]{geometry}               
\usepackage[parfill]{parskip}
\usepackage{xfrac}
\usepackage[font=small]{caption}
\usepackage{amsfonts}
\usepackage[intlimits]{amsmath}

\usepackage[section]{placeins}

\usepackage{mdwlist}
\usepackage{array}
\usepackage{multirow}
\usepackage{setspace}
\usepackage{graphicx}

\usepackage{float}
\restylefloat{figure}
\usepackage{sidecap}

\usepackage{tikz}
\usetikzlibrary{shapes,arrows}
\tikzset{myarrow/.style={draw,fill=black,double arrow,minimum height = 3.5cm,} } 
\tikzstyle{every picture}+=[remember picture]

\everymath{\displaystyle}

\title{A Generalization of a Bin-Based Modification \\ to the Stochastic Simulation Algorithm}
\author{David Collins }%\\ davidc@uvic.ca }
\affil{Dept. of Mathematics and Statistics, University of Victoria \\ PO Box 1700 STN CSC, Victoria, B.C. Canada V8W 2Y2}
\begin{document}
\maketitle

\textbf{Abstact} \\
In a 1996 paper, See$\beta$elberg, Trautmann and Thorn \cite{MS96} modified Gillespie's \cite{DG75b} Monte Carlo algorithm which is used to stochastically simulate the collision and coalescence process.  Their modification reduces the storage requirements of the simulation by several orders of magnitude.  However, their modification creates unphysical and potentially fatal conditions when used with common initial distributions.  We identify those conditions and propose a solution to maintain physically real conditions for all state variables during the evolution of any initial distribution.


%textbf{Introduction} 
\section{Introduction}
Gillespie \cite{DG75b} presented a Monte Carlo algorithm that stochastically simulates the collision and coalescence of liquid droplets.  His algorithm required the recording of information about $N$ droplets at each time step which demands impractically large storage requirements.  See$\beta$elberg et. al. \cite{MS96} reduced the required storage by discretizing the droplet size spectrum such that $N$ droplets are contained in $n$ bins.  Each bin contains of an aggregate of droplets.  Thus, bin $i$ contains an aggregate mass $m_i$ and a aggregate droplet number $N_i$.  The modification of See$\beta$elberg et. al. requires the storage of $2n$ variables (mass and number for each bin) rather than the masses of $N$ droplets.  Storage requirements are significantly reduced since $N \gg 2n$ by $\sim O(6)$.
   
Gillespie accounted for every drop in his simulation.  This demanded that droplet numbers be non-negative integers.  He described his simulation as encompassing an entire cloud which suggests the interpretation of the droplet mass as being represented by the `continuous model' of coalescence where the evolution of the droplet distribution tracks the mass of every droplet \cite{DG75a}.  His exact method was theoretical.  Although it did provide algorithms to execute the theory, the computational requirements of the algorithms were not quantifiably scaled to actual clouds as was his theory.  The computational requirements of his exact method when scaled to simulate entire clouds cannot be met by the most powerful computers of today.

However, a stochastic simulation can model droplet collisions in a subset of a cloud, and the results of such a simulation could then be scaled over the entire cloud.  The simulated volume could be considered a spatial average.  Assuming ergodicity as shown by Pope \cite{JW10} an ensemble of simulations, and subsequent spatial average, suggested by Gillespie \cite{DG75a} would satisfy the expected variation over a homogeneous cloud, but would not account for entrainment or detrainment.  Such an interpretation would be consistent with `quasi-stochastic' or `pure stochastic' models of coalescence where the spectral evolution tracks the number of drops of each mass or the probability of any number of drops of each mass, respectively, as explained in one of Gillespie's 1975 papers \cite{DG75a}.

See$\beta$elberg and collaborators apparently followed Gillespie's precedent and developed their bin-based modification using the `continuous model' of coalescence which assigned non-negative integers to  aggregate droplet numbers $N_i$ for each bin.  They did not specify that aggregate droplet numbers $N_i$ in bin $i$ would be restricted to integer values.  However, non-integer values can create unphysical droplet values for bins with $1<N_i<2$ as is shown in Section \ref{sec:Existence}.  A simulation volume containing non-integer values of number concentrations would be consistent with `quasi-stochastic' or `pure stochastic' models of coalescence \cite{DG75a}.  They used an initial gamma distribution.  To construct such a distribution using integer values for droplet numbers requires a crude approximation of the gamma distribution.  An approximation of, otherwise smooth, initial distributions may create bias in the evolved distribution due to non-uniform adjustments to intra-bin number concentrations.  Whether they made such crude approximations for initial distributions or ran enough simulations that a reasonably sized subset was successful is unclear.  Only the latter option creates non-physical state spaces as described in Section \ref{sec:Existence}.

The Monte Carlo algorithm and the subsequent bin-based modification are methods used to simulate an exact evolution of moments of a probability distribution that may otherwise be evolved by the kinetic collection equation.  A review of the literature identifies the See$\beta$elberg et. al. 1996 paper as a reference when validating the use of bin-based moment methods such as Bott's Linear Flux Method \cite{MP02, AK00}, when identifying a Monte Carlo method \cite{AK04, AB98,SS09}, or when quantifying a coalescence efficiency \cite{MP98, SS09}.  In 1999 Trautmann, a coauthor of the 1996 article, expanded upon that paper, with other collaborators, to include an exponential in-bin distribution and log-radius coordinates \cite{TT99}.  In 1999 Tzivion and collaborators critiqued the work of See$\beta$elberg et. al. by noting that it did not include sensitivity tests w.r.t. number of bins, kernel used, the time step, or initial conditions \cite{ST99}.   Other than the 1999 article by Trautmann et. al., a reproduction of the bin-based modification to the stochastic simulation has not been found in the literature by this author.

Generalizing the bin-based modification of Gillespie's stochastic simulation to include positive non-integer droplet numbers is consistent with a `quasi-stochastic' or `pure stochastic' model  of coalescence  \cite{DG75a}.  When a stochastic algorithm uses non-integers to quantify droplet numbers in simulations of a subset of a physical cloud, the evolved droplet size information can be interpreted as the number of drops of each mass or the probability of any number of drops of each mass.  In either case by the ergodicity shown by Pope \cite{JW10}, an ensemble of simulations of the generalized bin-based, modified stochastic algorithm gives a mean and variance of droplet masses that can be interpreted to represent a homogeneous portion of a cloud.  

These first and second moment statistics can be used to verify the results of detailed spectral methods.  A discrete spectral method is, in turn, used to validate bulk parameterizations \cite{MK00,AS01}.  Accurate bulk parameterizations are used by cloud resolving models to determine the onset of precipitation and the amount of radar reflectivity \cite{AS01,CF08}.  To make this necessary process more widely accessible we present a further modification that eliminates conditions leading to unphysical state spaces.

The discretization process introduces physical requirements on the modified evolution algorithm that were not present when individual droplets were selected by Gillespie's Monte Carlo Method.  The mean mass within a bin is defined (for bin $i$) by $\overline x_i = \sfrac{m_i}{N_i}$.  For the rest of this paper the phrase `mean mass' will refer to the mean mass within a bin.  Each collision in a simulation of the modified Gillespie algorithm involves reducing the droplet number of each selected source bin by unity and choosing, for removal, a droplet mass within the source bin.  The droplet mass must be chosen so that the remaining mass and remaining number yields a mean mass that is within the bin boundaries.  In Section \ref{sec:CurrentSelection}  the droplet selection process of the modified Monte Carlo algorithm by See$\beta$elberg and collaborators is briefly reviewed.  In Section \ref{sec:UnphysicalConditions} the existence and prevalence of conditions leading to unphysical state spaces are presented ,and an example is given.  Section \ref{sec:Refinement} details the modification to the droplet selection process used to eliminate unphysical conditions. 

%{Current Droplet Selection Process}
\section{Current Intra-Bin Droplet Selection Process}  \label{sec:CurrentSelection}
The algorithm by See$\beta$elberg et. al. includes a process to select a subinterval within the bin from which to choose a droplet mass.  The subinterval is given in their paper as Equation (9).  See$\beta$elberg et. al. assumed a uniform distribution of mass density within this subinterval and randomly selected a droplet.  In Sections \ref{sec:UnphysicalConditions} and \ref{sec:Refinement} we maintain that assumption.  The interval $I_i$ is constructed to be symmetric about the mean mass $\overline{x}_i$.  The masses of droplets at the left and right-hand bin boundaries are denoted as $x_{i-1/2}$ and $x_{i+1/2}$, respectively.
Letting $\Delta m_i \equiv \min(x_{i+1/2} - \overline x_i , \overline x_i - x_{i-1/2} )$, the bounds of $I_i$ are defined by
\begin{center}	
$d_{i-1/2}=\overline x_i - \Delta m_i, $, \hspace{1.0cm} and \hspace{1.0cm} $d_{i+1/2} = \Delta m_i + \overline x_i$ \hspace{1.0cm} \cite{MS96}
\end{center}
\begin{equation}
I_i = [d_{i-1/2} , d_{i+1/2} ] \subseteq [x_{i-1/2},x_{i+1/2}]  \nonumber
\end{equation}
 
\section{Unphysical Conditions} \label{sec:UnphysicalConditions}
The existence if the conditions that facilitate the unphysical state space is discussed.  An example which arises from the initial distribution used by See$\beta$elberg et. al. is detailed, and the prevalence of the condition is noted.

%{Unphysical Conditions: Existence} 
\subsection{Existence} \label{sec:Existence}
The droplet selection process given by See$\beta$elberg et. al. fails to adequately restrict the interval $I_i$ when $1 < N_i < 2$.  An example of the conditions predicating this problem are present in the initial droplet size distribution used by See$\beta$elberg et. al., and commonly used by others \cite{EB67,EB74,AB98,CF08}.  This initial gamma distribution is given in the paper by See$\beta$elberg et. al. as their Equation (19).  When for example bin $i$ containing $N_i \in (1,2)$ droplets is selected as the bin containing a pre-collision droplet for the next collision, the removal of too much or too little mass is possible and would result in a post-collision mean mass for bin $i$ that is either less than the mass of the smaller bin boundary, or greater than mass of the larger bin boundary, respectively.

\begin{figure} [htb] %
\centering{
\begin{tikzpicture}
\draw [<->,ultra thick] (0,0) -- (8,0);
\draw [-] (1,-0.25) -- (1,0.25); 				% left bin edge
\node [below] at (1,-0.5) {$x_{i-1/2}$};
\draw [red] (6.85,0.5) arc (30:-30:1cm);						% left edge Seesselberg interval
\draw [red] (3.15,-0.5) arc (210:150:1cm);						% right edge Seesselberg interval
\draw [-] (5,-0.25) -- (5,0.25); 				% mean mass
\node [below] at (5,-0.4) {$\overline x_{i}$};
\draw [-] (6,-0.25) -- (6,0.25); 				% mass in bin
\node [below] at (6,-0.4) {$\hat g_{i}$};
\draw [-] (7,-0.25) -- (7,0.25); 				% right bin edge
\node [below] at (7,-0.5) {$x_{i+1/2}$};
\end{tikzpicture} }
\caption{A schematic of one bin on the droplet size spectrum showing the bin boundaries for bin $i$, mean intra-bin mass $\overline x_i$, and the aggregate mass contained in bin $i$, $\hat g_i$.  The larger (red) arcs are the boundaries of the interval within which a droplet is chosen for the next collision.  The algorithm to determine these boundaries are by See$\beta$elberg et. al. (1996) and briefly explained in Section \ref{sec:CurrentSelection}. \protect\label{fig:GS} } 
\end{figure}
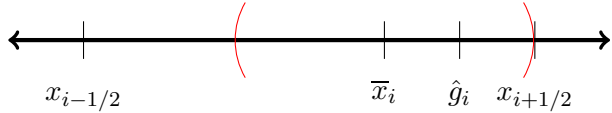

The sketch in Figure ($\ref{fig:GS}$) illustrates conditions that can lead to a mean mass outside of the bin after the collision.  The spectrum is described by the mass of a droplet where $\hat g_{i}$ is the total mass in bin $i$.  The two larger (red) arcs are the bounds of interval $I_i$.  A sub-interval bounded by $\hat g_i$ and $x_{i+1/2}$ lies within the interval bounded by the larger (red) arcs.  Therefore, when $N_i$ is slightly greater than unity, say $N_i=1+\alpha$ it is, (i) possible to select an amount of mass to be removed that is greater than the total amount of mass in bin $i$, and (ii) possible to select an amount of mass that is less than $\hat g_i$ but will result in the post-collision mean mass less than $x_{i-1/2}$ as further detailed in the following example. 

%{Unphysical Conditions: Example} 
\subsection{Example} \label{sec:Example}
The initial conditions presented in See$\beta$elberg et. al. create conditions in bin $\#17$ of their spectral discretization scheme that potentially lead to a post collision mean mass outside of the bin boundaries even when the mass removed is restricted to an amount less than the mass present.  The scaling of the initial distribution by 1500 as reported in their article only serves to shift the bin with $1<N_j<2$ downstream.  Bin $\#17$ of their unscaled distribution contains $2.046\times10^{-08}$g and 1.089 droplets.  When $1.893\times10^{-08}$g are removed then the remaining mass $x_{17,r}$ and remaining droplet number 0.089 gives a mean mass equal to the mass of the smaller bin boundary.  Therefore, when a droplet of mass $x_{17,\text{r}}$ is removed from bin $\#17$ and $x_{17,\text{r}} > 1.893\times10^{-08}$g, the post collision mean mass will be greater than the upper bound of the bin $\#17$.  

This flaw can be expressed as $\hat x_{i}< x_{i,\text{r}} <\hat g_{i}$ where $\hat x_{i}$ is the maximum mass that can be removed while keeping the post collision mean mass inside of the bin.  Having a mean mass outside of a bin that contains a fraction of a droplet such as 0.089 may not seem critical, but the addition of a droplet to bin $i$, due to a collision with a smaller droplet, will not guarantee that the mean mass will return to the inside of the bin.  

%{Unphysical Conditions: Prevalence} 
\subsection{Prevalence} \label{sec:Prevalence}
The prevalent use of a gamma function as an initial distribution requires that bins with $1 < N_i < 2$ will typically be present.  The persistence of a practical maximum droplet size (whether through precipitation or breakup) guarantees the existence of the condition in some bin at each time step.  A distribution which develops a bimodal character can have multiple bin with potential fatalities.  The stochastic nature of the algorithm facilitates a random selection of a source bin, and the possibility that a bin with a critically low droplet number be chosen.

Several simulations have been run without the refinement presented in this paper, and many have produced a mean mass external to the bin.  These external mean mass occurrences have occurred in multiple bins and at various time steps.  Successfully completing a simulation without a fatal error is possible.  However, such results may be produced after the existence of an unphysical state space at some time step during the simulation which may call into question the validity of the simulation.

%\section{A Refinement to the Droplet Selection Process}  
\section{A Refinement to the Droplet Selection Process}  \label{sec:Refinement}
\begin{figure} [htb]
\centering{
\begin{tikzpicture}
\draw [<->,ultra thick] (0,0) -- (8,0);
\draw [-] (1,-0.25) -- (1,0.25); 				% left bin edge
\node [below] at (1,-0.5) {$x_{i-1/2}$};
\draw [red] (6.85,0.5) arc (30:-30:1cm);						% left edge Seesselberg interval
\draw [red] (3.15,-0.5) arc (210:150:1cm);						% right edge Seesselberg interval
\draw [-] (5,-0.25) -- (5,0.25); 				% mean mass
\node [below] at (5,-0.4) {$\overline x_{i}$};
\draw [-] (6,-0.25) -- (6,0.25); 	
\draw [blue] (5.75,0.25) arc (15:-15:1cm);	
\node [above] at (5.75,0.4) {$\hat x_{i}$};
\node [below] at (6,-0.4) {$\hat g_{i}$};		% mass in bin
\draw [-] (7,-0.25) -- (7,0.25); 				% right bin edge
\node [below] at (7,-0.5) {$x_{i+1/2}$};
\end{tikzpicture} }
\caption{A schematic of one bin on the droplet size spectrum showing the same as in Figure ($\ref{fig:GS}$) with the addition of a smaller (blue) arc $\hat x_i$ that represents the upper bound of a more restrict interval.  The calculations used to create this more restrictive boundary are described in the surrounding text. \protect\label{fig:GS2} }
\end{figure}
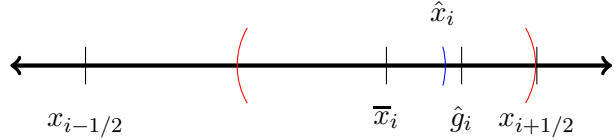

A sub-interval $\hat I_i$ is constructed that is dependent on the bin boundaries ($x_{i-1/2}$ and $x_{i+1/2}$): 
$\hat I_i \subseteq I_i$ such that after the simulation of a collision is complete the mean mass $\overline x_i$ will remain between the bin boundaries.  The upper bound of $\hat I_i$ is shown in Figure ($\ref{fig:GS2}$) and subsequently developed herein.

The smaller (blue) arc represents the maximum amount of mass that can be removed such that $\hat g_{i} - \hat x_{i}$ is the amount of mass remaining after the next collision involving bin $i$ and $\left({\hat g_{i} - \hat x_{i}}\right)/{\alpha}=\overline x_{i,\text{a}}=x_{i-1/2}$ where $\overline x_{i,\text{a}}$ is the mean mass after the collision.  The following procedure is used to calculate $\hat x_{i}$ as a function of $x_{i-1/2}$ and also to calculate a more restrictive lower bound as a function of $x_{i+1/2}$.

The following procedure will give one or more bounds of the more restrictive interval that ensures the remaining mean intra-bin mass will be within the bin boundaries:
\begin{equation}
e_{i-1/2} = L_{i}-x_{i+1/2}(N_{i}-1)  \hspace{1.5cm}  e_{i+1/2} = \hat x_i = L_{i}-x_{i-1/2}(N_{i}-1)  \nonumber
\end{equation}
\begin{equation}
\overline d_{i-1/2} = \max(d_{i-1/2},e_{i-1/2}), \hspace{1.5cm}   \overline d_{i+1/2} = \min(d_{i+1/2},e_{i+1/2}) \nonumber
\end{equation}
\begin{equation}
\overline I_i = [\overline d_{i-1/2},\overline d_{i+1/2}] \subseteq I_i \subseteq [x_{i-1/2},x_{i+1/2}] 
\nonumber
\end{equation}
$L_{i}$ and $N_{i}$ are the liquid water content and number concentration in bin $i$ prior to the next collision involving bin $i$.  

In the example from Section \ref{sec:Example}, $\overline d_{i+1/2}$ is a function of the lower bin boundary and is shown in Figure (\ref{fig:GS3}).

Notice that the interval bounded by the left-most of the larger (red) arcs and $\hat x_{i}$ is not symmetric about $\overline x_{i}$.  An interval symmetric about $\overline x_{i}$ is needed to avoid bias.  Therefore, a final step in the spirit of the algorithm by See$\beta$selberg et. al. (1996) is used to produce the sub-interval $\hat I_i$ bounded by the two smaller (blue) arcs:
\begin{equation}
\hat x'_{i} = \overline x_{i} - (\hat x_{i} - \overline x_{i})  \nonumber
\end{equation}
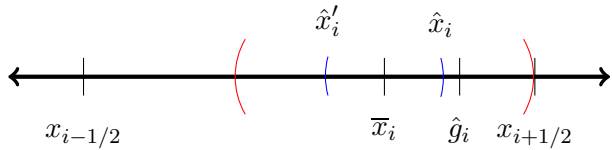
\begin{figure} [htb]
\centering{
\begin{tikzpicture}
\draw [<->,ultra thick] (0,0) -- (8,0);
\draw [-] (1,-0.25) -- (1,0.25); 				% left bin edge
\node [below] at (1,-0.5) {$x_{i-1/2}$};
\draw [red] (3.15,-0.5) arc (210:150:1cm);						% left edge Seesselberg interval
\draw [blue] (4.25,-0.25) arc (195:165:1cm);	
\node [above] at (4.25,0.4) {$\hat x'_{i}$};
\draw [-] (5,-0.25) -- (5,0.25); 				% mean mass
\node [below] at (5,-0.4) {$\overline x_{i}$};
\draw [-] (6,-0.25) -- (6,0.25); 	
\draw [blue] (5.75,0.25) arc (15:-15:1cm);	
\node [above] at (5.75,0.4) {$\hat x_{i}$};
\node [below] at (6,-0.4) {$\hat g_{i}$};		% mass in bin
\draw [red] (6.85,0.5) arc (30:-30:1cm);						% right edge Seesselberg interval
\draw [-] (7,-0.25) -- (7,0.25); 				% right bin edge
\node [below] at (7,-0.5) {$x_{i+1/2}$};
\end{tikzpicture} }
\caption{A schematic of one bin on the droplet size spectrum showing the same as in Figure ($\ref{fig:GS2}$) with the addition of a second smaller (blue) arc $\hat x'_i$ that represents the lower bound of a more the restrict interval.  The calculations used to create the lower bound of this more restrictive interval are in the same spirit as that originally used by See$\beta$selberg et. al. (1996).} \label{fig:GS3}
\end{figure}
We now have an interval $\hat I_i=[\hat x'_{i},\hat x_{i}]$, symmetric about $ \overline x_{i}$, from which to draw the pre-collision droplet which guarantees that the remaining mass and remaining number concentration will yield a post-collision mean intra-bin mass that is inside of the bin boundaries:
\begin{center}
$\hat I_i = [\hat x'_{i},\hat x_{i}] \subseteq \overline I_i \subseteq I_i \subseteq [x_{i-1/2},x_{i+1/2}]$
\end{center}

\section{Conclusion}  \label{sec:Conclusion}
The discretization of the droplet spectrum by See$\beta$elberg (1996) and collaborators reduced the computational expense of stochastically simulating the collision and coalescence process over the expense of a stochastic simulation that records information about every droplet as given by Gillespie (1975).  That modification produce the additional requirement of ensuring that the ratio of intra-bin mass and intra-bin number lie between the two values on the edges of the bin.  This requirement is intrinsically met by the algorithm presented by See$\beta$elberg et. al. except for the cases when the post-collision intra-bin number is $0< N_i<1$.  They restrict the interval in the source bins for the purpose of eliminating droplet selection bias.  However, their selection interval does not guarantee that the post-collision mean mass lies within the bin.  Situations resulting in such a post collision number concentration has been shown to occur in the commonly used gamma initial distribution, and the condition will persist as the mode traverses the droplet spectrum. 

These situations are addressed by further restricting the intra-bin interval contained in the source bins upon which is constructed a uniform distribution.  The additional restrictions on the domain from which the mass of a source droplet can be selected are functions of bin boundaries.  The new domain is then made symmetric about the intra-bin mean mass to eliminate selection bias.  The restricted intra-bin domain allows for a generalization of initial conditions that include continuous  (i.e. non-integer) values of number concentration.  There is no need to modify an initial distribution that contains non-integer values of number concentration.  Nor is there a need to prevent collisions of droplets in bins containing a number between 1 and 2. \\

%%%%%
\bibliographystyle{plain}
\bibliography{GSSbib}

\end{document}